\documentclass[12pt]{amsart}
\usepackage{amsmath,amssymb,amsbsy,amsfonts,latexsym,amsopn,amstext,cite,
                                               amsxtra,euscript,amscd,bm}
\usepackage{url}
\usepackage{mathrsfs}

\usepackage{color}
\usepackage[colorlinks,linkcolor=blue,anchorcolor=blue,citecolor=blue,backref=page]{hyperref}
\usepackage{color}
\usepackage{graphics,epsfig}
\usepackage{graphicx}
\usepackage{float}
\usepackage{epstopdf}
\hypersetup{breaklinks=true}

\usepackage[top=1in, bottom=1in, left=1in, right=1in]{geometry}

\usepackage[np]{numprint}
\npdecimalsign{\ensuremath{.}}

\def\Th1{\varTheta}

\usepackage{bibentry}
\usepackage{bbm}

\usepackage[english]{babel}
\usepackage{mathtools}
\usepackage{todonotes}
\usepackage{url}
\usepackage[colorlinks,linkcolor=blue,anchorcolor=blue,citecolor=blue,backref=page]{hyperref}

\usepackage[norefs,nocites]{refcheck}

\begin{document}

\newtheorem{theorem}{Theorem}
\newtheorem{lemma}[theorem]{Lemma}
\newtheorem{claim}[theorem]{Claim}
\newtheorem{cor}[theorem]{Corollary}
\newtheorem{conj}[theorem]{Conjecture}
\newtheorem{prop}[theorem]{Proposition}
\newtheorem{definition}[theorem]{Definition}
\newtheorem{question}[theorem]{Question}
\newtheorem{example}[theorem]{Example}
\newcommand{\hh}{{{\mathrm h}}}
\newtheorem{remark}[theorem]{Remark}

\numberwithin{equation}{section}
\numberwithin{theorem}{section}
\numberwithin{table}{section}
\numberwithin{figure}{section}

\def\sssum{\mathop{\sum\!\sum\!\sum}}
\def\ssum{\mathop{\sum\ldots \sum}}
\def\iint{\mathop{\int\ldots \int}}

\newcommand{\diam}{\operatorname{diam}}

\def\squareforqed{\hbox{\rlap{$\sqcap$}$\sqcup$}}
\def\qed{\ifmmode\squareforqed\else{\unskip\nobreak\hfil
\penalty50\hskip1em \nobreak\hfil\squareforqed
\parfillskip=0pt\finalhyphendemerits=0\endgraf}\fi}

\newfont{\teneufm}{eufm10}
\newfont{\seveneufm}{eufm7}
\newfont{\fiveeufm}{eufm5}
%
%
\newfam\eufmfam
     \textfont\eufmfam=\teneufm
\scriptfont\eufmfam=\seveneufm
     \scriptscriptfont\eufmfam=\fiveeufm
%
%
\def\frak#1{{\fam\eufmfam\relax#1}}

\newcommand{\bflambda}{{\boldsymbol{\lambda}}}
\newcommand{\bfmu}{{\boldsymbol{\mu}}}
\newcommand{\bfxi}{{\boldsymbol{\eta}}}
\newcommand{\bfrho}{{\boldsymbol{\rho}}}

\def\eps{\varepsilon}

\def\fK{\mathfrak K}
\def\fT{\mathfrak{T}}
\def\fL{\mathfrak L}
\def\fR{\mathfrak R}

\def\fA{{\mathfrak A}}
\def\fB{{\mathfrak B}}
\def\fC{{\mathfrak C}}
\def\fM{{\mathfrak M}}
\def\fS{{\mathfrak  S}}
\def\fU{{\mathfrak U}}
\def\fW{{\mathfrak W}}

\def\T {\mathsf {T}}
\def\Tor{\mathsf{T}_d}
\def\Tore{\widetilde{\mathrm{T}}_{d} }

\def\sM {\mathsf {M}}

\def\ss{\mathsf {s}}

\def\Kmnd{\cK_d(m,n)}
\def\Kmnp{\cK_p(m,n)}
\def\Kmnq{\cK_q(m,n)}

\def \balpha{\bm{\alpha}}
\def \bbeta{\bm{\beta}}
\def \bgamma{\bm{\gamma}}
\def \bdelta{\bm{\delta}}
\def \bzeta{\bm{\zeta}}
\def \blambda{\bm{\lambda}}
\def \bchi{\bm{\chi}}
\def \bphi{\bm{\varphi}}
\def \bpsi{\bm{\psi}}
\def \bnu{\bm{\nu}}
\def \bomega{\bm{\omega}}

\def \bell{\bm{\ell}}

\def\eqref#1{(\ref{#1})}

\def\vec#1{\mathbf{#1}}

\newcommand{\abs}[1]{\left| #1 \right|}

\def\Zq{\mathbb{Z}_q}
\def\Zqx{\mathbb{Z}_q^*}
\def\Zd{\mathbb{Z}_d}
\def\Zdx{\mathbb{Z}_d^*}
\def\Zf{\mathbb{Z}_f}
\def\Zfx{\mathbb{Z}_f^*}
\def\Zp{\mathbb{Z}_p}
\def\Zpx{\mathbb{Z}_p^*}
\def\cM{\mathcal M}
\def\cE{\mathcal E}
\def\cH{\mathcal H}

\def\le{\leqslant}

\def\ge{\geqslant}

\def\sfB{\mathsf {B}}
\def\sfC{\mathsf {C}}
\def\L{\mathsf {L}}
\def\FF{\mathsf {F}}

\def\sE {\mathscr{E}}
\def\sS {\mathscr{S}}

\def\cA{{\mathcal A}}
\def\cB{{\mathcal B}}
\def\cC{{\mathcal C}}
\def\cD{{\mathcal D}}
\def\cE{{\mathcal E}}
\def\cF{{\mathcal F}}
\def\cG{{\mathcal G}}
\def\cH{{\mathcal H}}
\def\cI{{\mathcal I}}
\def\cJ{{\mathcal J}}
\def\cK{{\mathcal K}}
\def\cL{{\mathcal L}}
\def\cM{{\mathcal M}}
\def\cN{{\mathcal N}}
\def\cO{{\mathcal O}}
\def\cP{{\mathcal P}}
\def\cQ{{\mathcal Q}}
\def\cR{{\mathcal R}}
\def\cS{{\mathcal S}}
\def\cT{{\mathcal T}}
\def\cU{{\mathcal U}}
\def\cV{{\mathcal V}}
\def\cW{{\mathcal W}}
\def\cX{{\mathcal X}}
\def\cY{{\mathcal Y}}
\def\cZ{{\mathcal Z}}
\newcommand{\rmod}[1]{\: \mbox{mod} \: #1}

\def\cg{{\mathcal g}}

\def\vy{\mathbf y}
\def\vr{\mathbf r}
\def\vx{\mathbf x}
\def\va{\mathbf a}
\def\vb{\mathbf b}
\def\vc{\mathbf c}
\def\ve{\mathbf e}
\def\vh{\mathbf h}
\def\vk{\mathbf k}
\def\vm{\mathbf m}
\def\vz{\mathbf z}
\def\vu{\mathbf u}
\def\vv{\mathbf v}

\def\e{{\mathbf{\,e}}}
\def\ep{{\mathbf{\,e}}_p}
\def\eq{{\mathbf{\,e}}_q}

\def\Tr{{\mathrm{Tr}}}
\def\Nm{{\mathrm{Nm}}}

 \def\SS{{\mathbf{S}}}

\def\lcm{{\mathrm{lcm}}}

 \def\0{{\mathbf{0}}}

\def\({\left(}
\def\){\right)}
\def\l|{\left|}
\def\r|{\right|}
\def\fl#1{\left\lfloor#1\right\rfloor}
\def\rf#1{\left\lceil#1\right\rceil}
\def\fl#1{\left\lfloor#1\right\rfloor}
\def\ni#1{\left\lfloor#1\right\rceil}
\def\sumstar#1{\mathop{\sum\vphantom|^{\!\!*}\,}_{#1}}

\def\mand{\qquad \mbox{and} \qquad}

\def\tblue#1{\begin{color}{blue}{{#1}}\end{color}}




\hyphenation{re-pub-lished}

\mathsurround=1pt

\def\bfdefault{b}

\def \E{{\mathbb E}}
\def \F{{\mathbb F}}
\def \K{{\mathbb K}}
\def \N{{\mathbb N}}
\def \Z{{\mathbb Z}}
\def \P{{\mathbb P}}
\def \Q{{\mathbb Q}}
\def \R{{\mathbb R}}
\def \C{{\mathbb C}}
\def\Fp{\F_p}
\def \fp{\Fp^*}

 \def \xbar{\overline x}

\title[Large Gaps PS Primes]
{Large Gaps Between Piatetski-Shapiro Primes}

\author [L. P. Wijaya]{Laurence P. Wijaya}
\address{Department of Mathematics, University of Kentucky, 715 Patterson Office Tower, Lexington, KY 40506, USA}
\email{laurence.wijaya@uky.edu}

\begin{abstract}  
A Piatetski-Shapiro prime is a prime number that can be written as $\lfloor m^c \rfloor$ for some positive integer $m$ and non-integer $c>1$. Li and Pan showed the results for small gaps between consecutive Piatetski-Shapiro primes for $1<c<9/8$. On the other hand, Deshouillers, Drmota, and M\"ullner showed the large gaps between Piatetski-Shapiro primes although they did not mention it explicitly whenever there are infinitely many Piatetski-Shapiro primes with exponent $c$ using binomial expansion. In this article, we use a different method from theirs, and we obtain the explicit exponent that appears in their paper.
 \end{abstract}


\maketitle


\section{Introduction}
One interesting subset of primes is the set of Piatetski-Shapiro primes, introduced by Piatetski-Shapiro in \cite{PS}, which is defined by
\[
\mathcal{P}_c=\{p\text{ prime }: p=\lfloor m^{(c)}\rfloor\text{ for some }m\in\N\}
\]
with $c>1$ is a positive real number that is not an integer. The case $1<c<2$ is in general quite well-studied. One of the reasons that this set is interesting is because this is a relatively sparse subset of primes with an explicit formula, and we can expect that many interesting properties for the set of all prime numbers also hold for the set of Piatetski-Shapiro primes.

One interesting known fact is that $\mathcal{P}_c$ is infinite in the range $1<c<243/205$ from \cite{RW}. In fact, Guo, Guo, and Lu \cite{GGL} gives the best range for $c$ for the prime number theorem for Piatetski-Shapiro primes, but the range of $c$ is smaller, which is $1<c<10318869/8886224$. We hope that we can get the results for all $1<c<2$ although these are still unknown.

Besides the question on the infinitude of Piatetski-Shapiro primes, or the exact asymptotic order of their distribution, one can ask about the gaps between consecutive Piatetski-Shapiro primes. Clearly because of how they are defined, the gaps between consecutive Piatetski-Shapiro primes cannot be bounded since we have that $(m+1)^c-m^c$ goes to infinity as $m$ goes to infinity for $c>1$, unlike the case of the set of all primes. See \cite{S} for the most recent update on bounded gaps between primes.

Li and Pan \cite{LP}, using modification of the Maynard-Tao sieve, showed that for $1<c<9/8$, if $\mathcal{P}_c=\{q^{(c)}_1,q^{(c)}_2,\dots\}$ with $q_1^{(c)}<q_2^{(c)}<\dots$, we have
\[
\liminf_{n\rightarrow\infty} \frac{q_{n+1}^{(c)}-q_n^{(c)}}{\left(q_n^{(c)}\right)^{1-1/c}}<C_c
\]
for some constant $C_c$ that depends only on $c$. Note that the average gaps between consecutive Piatetski-Shapiro primes up to $x$, when we have the prime number theorem for Piatetski-Shapiro primes, is roughly
\[
x^{1-1/c}\log x,
\]
so this means that the limit infimum of the gaps of consecutive Piatetski-Shapiro primes is smaller than the average. One can view this as similar to the result for the whole primes, where we do not have the logarithmic term.

On the other hand, we can ask about the other gap, which is the limit supremum of the consecutive gaps. In this case, it turns out that we have stronger result compared to the whole primes. For the set of all primes, Ford, Green, Maynard, and Tao \cite{FGMT} showed that
\[
\max_{p_n\leq X} (p_{n+1}-p_n)\gg (c+o(1))\frac{\log X \log \log X\log\log\log\log X}{\log\log\log X}
\]
for some positive constant $c$, where $p_1<p_2<\dots$ is the enumeration of all prime numbers in increasing order.

However, in the case of Piatetski-Shapiro primes, we have the following larger gap.
\begin{theorem}
\label{thm : main}
    There are infinitely many $n$ such that for $1<c<243/205$
    \[
    q_{n+1}^{(c)}-q_n^{(c)}\gg_c\sqrt{q_n^{(c)}}.
    \]
\end{theorem}

Note that a similar result already appeared in Deshouillers, Drmota, and M\"ullner's paper \cite{DDM}. We refer to Remark \ref{rem:DDM} at the end of this article for more detail.

We remark that one of the conjectures for prime numbers, called Legendre's conjecture that says there is at least one prime between $m^2$ and $(m+1)^2$ for any $m>1$, would imply that
\[
p_{n+1}-p_n<\sqrt{p_n} 
\]
asymptotically. See \cite{CS,DJ} for related topics. Thus, for Piatetski-Shapiro primes, we do not have such result even if the Legendre's conjecture is true.

The outline of the proof is as follows. Instead of using sieve method, we use a lemma stating that we can write the tangent line of $f(x)=x^c$ such that its constant term that differs slightly from an even integer, and the slope is even integer at large $x$. Consequently, we can construct blocks of intervals such that every integer of the form $\lfloor  n^c \rfloor$ is even on each interval, and the length of the interval is as large as we want in the statement of Theorem \ref{thm : main}.

\section{Notations and Conventions}
We denote the set of all positive integers by $\N$. For two functions $f$ and $g$, we write $f\ll g$ or $f=O(g)$ or $g\gg f$ if
\[
|f(x)|\leq Cg(x)
\]
as $x\rightarrow\infty$ for some positive constant $C$. If the constant depends on another parameter, for example $c$
, we denote it by $f\ll_c g$ or $f=O_c(g)$ or $g\gg_c f$ and the constant that depends on $c$ will have a subscript, for example $C_c$.

For a sequence of real numbers $\{a_n\}$, by $\{a_n\mod 1\}$ we mean the sequence of fractional parts of $\{a_n\}$, so $a_n\mod 1\in [0,1)$ for any $n$.

\section{Proof of Theorem \ref{thm : main}}
To prove Theorem \ref{thm : main}, we need the following lemma that is inspired by \cite{SY} where they work on arithmetic progressions on the set $\{\lfloor n^c \rfloor : n\in \N\}$, which is bigger than the set of Piatetski-Shapiro primes. One can view this as a lemma about the tangent line of $f(x)=x^c$.
\begin{lemma}
\label{lem:tangent}
    Let $f(x)=x^c$. Then for any $T>0$, there exists arbitrary large real number $t>T$ such that
    \[
    f(t)+f'(t)(x-t)=Ax+B
    \]
    where $A$ is an even integer, $B=2\ell+\theta$, $\ell\in\Z, 0\leq\theta\leq 1/2$, and for any $x$ satisfying $|x-t|\leq H$ with $H=o(t)$.
\end{lemma}
\begin{proof}
    We divide this into two cases, depending on whether $c$ is rational or irrational.

    \textbf{Case 1.} $c$ is rational.

    Let $c=a/b$ with $\gcd(a,b)=1$. Choose a large positive integer $v$ such that it is divisible by $2b$, and set $t=v^b$ such that it is larger than the given $T$. Now if we define
    \[
    A:=\frac{a}{b}v^{a-b},\quad B=f(t)-At,
    \]
    we see that
    \[
    f(t)+f'(t)(x-t)=v^a+\frac{a}{b}v^{a-b}(x-v^b)=Ax+B
    \]
    where in this case, $A,B$ both are even integer since we have that $v$ is an even integer.

    \textbf{Case 2.} $c$ is irrational.

    In this case, for $k\in\N$, define $t_k=(2k/c)^{1/(c-1)}$. We will choose suitable $k$ later. We claim that the assertion holds with $A=2k$ and $B=f(t_k)-2kt_k$. We can choose $k$ large such that $t_k>T$ for the given $T$. Since $f'(t_k)=2k$, our claim holds for the value of $A$. Now we have
    \[
    B=-(c-1)\left(\frac{2k}{c}\right)^{c/(c-1)}.
    \]
    We can write
    \[
    B/2=-\left(\frac{c-1}{2}\right)\left(\frac{2}{c}\right)^{c/(c-1)}k^{c/(c-1)}.
    \]
    Using equidistribution of the sequence $\left\{\left(\frac{c-1}{2}\right)\left(\frac{2}{c}\right)^{c/(c-1)}k^{c/(c-1)}\mod 1\right\}_{k\geq 1}$ that is known from \cite{C} (see also \cite{TY} for some history on the topic of equidistribution of the sequence of this type), there are infinitely many $k$ such that
    \[
    \left\{-\left(\frac{c-1}{2}\right)\left(\frac{2}{c}\right)^{c/(c-1)}k^{c/(c-1)}\mod 1\right\}\in [1/8,1/4].
    \]
    Thus, we can choose large enough $k$ such that $t_k>T$, $-\left(\frac{c-1}{2}\right)\left(\frac{2}{c}\right)^{c/(c-1)}k^{c/(c-1)}\mod 1\in [1/8,1/4]$. Taking our $t=t_k$ then gives $B=2\ell+\theta$ with $1/4\leq\theta\leq 1/2$.
\end{proof}
Now we consider the Taylor expansion of $f$. We need to control the error term after the first derivative to get our main result.
\begin{lemma}
\label{lem : Taylor}
    Fix a large $t$ that satisfies the previous lemma. Uniformly for $|x-t|\leq H$ with $H=o(t)$, we have
    \[
    f(x)=Ax+B+E_t(x),\quad 0\leq E_t(x)\leq K_c t^{c-2}H^2.
    \]
\end{lemma}
\begin{proof}
    We use Taylor's expansion around $t$ and Lemma \ref{lem:tangent}. We have since $f$ is strictly convex and strictly increasing for $x>0$,
    \[
    f(x)=f(t)+f'(t)(x-t)+E_t(x)=Ax+B+E_t(x)
    \]
    where $0\leq E_t(x)\leq K_c t^{c-2}H^2$ with $H=o(t)$ and $|x-t|\leq H.$
\end{proof}

Now we are ready to prove our main theorem.

\begin{proof}[Proof of Theorem \ref{thm : main}]
For any integer $|m-t|\leq H$ with $H=o(t)$, we have
\[
\lfloor m^c \rfloor=Am+2\ell+\lfloor \theta+E_t(m) \rfloor=Am+2\ell.
\]
So, if we let $a=\lceil t-H \rceil, b=\lfloor t+H\rfloor$ under the same assumptions as in Lemma \ref{lem : Taylor}, there is no Piatetski-Shapiro prime in $[\lfloor a^c \rfloor,\lfloor b^c \rfloor]$. 

Choose $H=\eta_ct^{1-c/2}$ with $\eta_c>0$ small such that $K_c\eta_c^2<1/4$ ensuring that $E_t(m)\leq 1/4$ in Lemma \ref{lem : Taylor}. So we have
\[
\lfloor b^c \rfloor-\lfloor a^c \rfloor\gg t^{c/2}.
\]
By choosing $q^{(c)}_n,q^{(c)}_{n+1}$ such that $[\lfloor a^c \rfloor,\lfloor b^c \rfloor]\subseteq [q^{(c)}_n,q^{(c)}_{n+1}]$, we obtain
\[
q_{n+1}^{(c)}-q_n^{(c)}\gg_c t^{c/2}\gg_c \sqrt{q_n^{(c)}}
\]
as claimed.
\end{proof}
\begin{remark}
    By a similar method as in the proof of Theorem \ref{thm : main}, we can take any positive integer $d>1$ and build blocks of intervals, with the length of each interval being $\gg t^{c/2}$ such that for every integer of the form $\lfloor n^c \rfloor$ in the interval is divisible by $d$.
\end{remark}
\begin{remark}
\label{rem:DDM}
    The fact that for any $c\in (1,\infty)\backslash \N$ there exist $\kappa>0$ and infinitely many $n\in \N$ such that for every integer $h\in [0,n^{\kappa}]$ we have the numbers $\lfloor (n+h)^{c} \rfloor$ are even positive integers was established by Deshouillers, Drmota, and M\"ullner in \cite{DDM}. Their method was to consider the binomial coefficient of $(n+h)^c$, then use the Erd\H{o}s-Tur\'{a}n discrepancy result to conclude. Because of this, they did not get the explicit exponent $\kappa$.
\end{remark}
\section*{Acknowledgement}
The author would like to thank to Fernando Xuancheng Shao for his advice and pointing out \cite{DDM}. The author would also like to thank to Muhammad Afifurrahman for the proofreading.
\section*{AI Disclosure}
no Anime is Injured by the author during the writing process.


\begin{thebibliography}{99}

\bibitem{CS} M. Chamberland, A. Straub, `Weakening the Legendre Conjecture', \textit{arXiv preprint arXiv:2602.22502}, (2026).

\bibitem{C} P., Csillag, `Über die gleichmässige Verteilung nicht ganzer positiver Potenzen mod. 1', Acta Litt. Sci. Szeged 5, 13–18 (1930).

\bibitem{DDM} J.-M. Deshouillers, M. Drmota and C. M\"ullner, `Coprimality of consecutive elements in a Piatetski-Shapiro sequence', in {\it Number theory in memory of Eduard Wirsing}, 91--98, Springer, Cham, (2023).

\bibitem{DJ} A.~W. Dudek and D.~R. Johnston, `Almost primes between all squares', J. Number Theory {\bf 278}, 726--745, (2026).

\bibitem{FGMT} K.~B. Ford et al., `Long gaps between primes', J. Amer. Math. Soc. {\bf 31}, no.~1, 65--105, (2018).

\bibitem{GGL} L. Guo, V. Z. Guo, and L. Lu, `The Piatetski-Shapiro prime number theorem', \textit{arXiv preprint arXiv:2505.10391}, (2025).

\bibitem{LP} H. Li and H. Pan, `The Green-Tao theorem for Piatetski-Shapiro primes', J. Funct. Anal. {\bf 291}, no.~1, Paper No. 111440, 25 pp., (2026).

\bibitem{PS} I.~I. Piatetski-Shapiro, `On the distribution of prime numbers in sequences of the form $[f(n)]$', Mat. Sbornik N.S. {\bf 33/75}, 559--566, (1953).

\bibitem{RW} J. Rivat and J. Wu, `Prime numbers of the form $[n^{(c)}]$', Glasg. Math. J. {\bf 43}, no.~2, 237--254, (2001).

\bibitem{SY} K. Saito and Y. Yoshida, `Arithmetic progressions in the graphs of slightly curved sequences', J. Integer Seq. {\bf 22}, no.~2, Art. 19.2.1, 25 pp, (2019).

\bibitem{S} J. Stadlmann, `Bounded gaps between primes', \textit{arXiv preprint arXiv:2608.31126}, (2026).

\bibitem{TY} N. Technau and N. Yesha, `On the correlations of $n^\alpha\bmod1$', J. Eur. Math. Soc. (JEMS) {\bf 25}, no.~10, 4123--4154, (2023).


\end{thebibliography}
\end{document}